

\baselineskip=14pt
\parskip=10pt

\font\eightrm=cmr8 
\font\eighttt=cmtt8
\magnification=\magstephalf
\def\C{{\cal C}}
\def\P{{\cal P}}

\def\1{{\overline{1}}}
\def\2{{\overline{2}}}
\parindent=0pt
\overfullrule=0in

\def\frac#1#2{{#1 \over #2}}
\bf
\centerline
{
Automatic Generation of Theorems and Proofs on Enumerating Consecutive-Wilf Classes
}

\rm
\bigskip
\centerline{ 
By {\it Andrew BAXTER}{$^1$}, {\it Brian NAKAMURA}{$^1$}, and
{\it Doron ZEILBERGER}\footnote{$^1$}
{\eightrm  \raggedright
Department of Mathematics, Rutgers University (New Brunswick),
Hill Center-Busch Campus, 110 Frelinghuysen Rd., Piscataway,
NJ 08854-8019, USA.
{\eighttt [baxter, bnaka, zeilberg]@math.rutgers.edu} ,
\hfill\break
{\eighttt http://www.math.rutgers.edu/[\~{}bnaka/, \~{}baxter/, \~{}zeilberg/]}  .
Jan. 20, 2011.
Accompanied by Maple packages {\eighttt ELIZALDE} and {\eighttt SERGI}
downloadable from the webpage of this article:
\hfill\break
{\eighttt http://www.math.rutgers.edu/\~{}zeilberg/mamarim/mamarimhtml/auto.html} \quad ,
\hfill \break
where  the reader can find lots of sample input and output.
The work of AB and DZ was supported in part by the United States of America National Science Foundation.
}
}

\quad \quad \quad\quad\quad\quad\quad\quad\quad\quad\quad\quad {\it To W from Z (et. al.), a gift for his $\frac{2}{3}|S_5|$-th birthday}

{\bf Preface}  

This article describes two complementary approaches to enumeration, the {\it positive}
and the {\it negative}, each with its advantages and disadvantages. Both approaches are amenable
to {\it automation}, and when applied to the currently active subarea, initiated in 2003 by
Sergi Elizalde and Marc Noy[EN], of {\it consecutive pattern-avoidance} in permutations, were successfully
pursued by the first two authors Andrew Baxter[B] and Brian Nakamura[N]. 
This article summarizes their research and in the case of [N] presents
an umbral viewpoint to the same approach.
The main purpose of this article is to
briefly explain the Maple packages, {\tt SERGI} and {\tt ELIZALDE},
developed by AB-DZ and BN-DZ respectively,
implementing the algorithms that enable the computer to ``do research" by
deriving, {\it all by itself}, functional equations for the generating functions
that enable polynomial-time enumeration for any set of patterns.
In the case of ELIZALDE (the ``negative" approach), these functional equations can be
sometimes (automatically!) simplified, and imply ``explicit" formulas, that
previously were derived by humans using ad-hoc methods. We also get lots of
new ``explicit" results, beyond the scope of humans, but we have to admit that we still
need humans to handle ``infinite families" of patterns, but this too, no doubt, will
soon be automatable, and we leave this as a challenge to the (human and/or computer) reader.

{\bf Consecutive Pattern Avoidance}

Inspired by the very active research in pattern-avoidance, pioneered by Herb Wilf, Rodica Simion,
Frank Schmidt, Richard Stanley, Don Knuth and others, Sergi Elizalde, in his PhD thesis (written under the
direction of Richard Stanley)
introduced the study of permutations avoiding {\it consecutive patterns}.

Recall that an $n$-{\it permutation} is a  sequence of integers $\pi=\pi_1 \dots \pi_n$ 
of length $n$ where each integer in $\{1, \dots , n\}$ appears exactly once.
It is well-known and very easy to see (today!) that the number of $n$-permutations is $n!:= \prod_{i=1}^{n} i \quad $.

The {\it reduction} of a list of different (integer or real) numbers (or members of any totally ordered set) $[i_1, i_2, \dots, i_k]$,
to be denoted by $R([i_1, i_2, \dots, i_k])$, is the permutation of $\{1,2, \dots, k \}$ that preserves the relative rankings of the entries.
In other words, $p_i<p_j$ iff $q_i<q_j$.
For example the reduction of $[4,2,7,5]$ is $[2,1,4,3]$ and the reduction of $[\pi, e, \gamma, \phi]$ is $[4,3,1,2]$.

Fixing a pattern $p=[p_1, \dots , p_k]$, a permutation $\pi=[\pi_1, \dots , \pi_n]$ {\it avoids} the consecutive pattern
$p$ if for  all $i$, $1 \leq i \leq n-k+1$, the reduction of the list $[\pi_i,\pi_{i+1}, \dots, \pi_{i+k-1}]$ is {\it not} $p$.
More generally a permutation $\pi$ avoids a set of patterns $\P$ if it avoids each and
every pattern $p \in \P$.

The central problem is to answer the question: ``Given a pattern or a set of patterns, find a `formula',
or at least an efficient algorithm (in the sense of Wilf[W]), that inputs a positive integer $n$ and outputs the number of
permutations of length $n$ that avoid that pattern (or set of patterns)''.

{\bf Human Research}

After the pioneering work of Elizalde and Noy [EN], quite a few people contributed significantly,
including Anders Claesson, Toufik Mansour,  Sergey Kitaev, Anthony Mendes, Jeff Remmel,
and more recently, Vladimir Dotsenko, Anton Khoroshkin and Boris Shapiro.
Also recently we witnessed the beautiful resolution of the Warlimont conjecture by
Richard Ehrenborg, Sergey Kitaev, and Peter Perry [EKP].
The latter paper also contains extensive references.

{\bf Recommended Reading}

While the present article tries to be self-contained, the readers would get more out of  it if
they are familiar with [Z1]. Other applications of the umbral transfer matrix method
were given in [EZ][Z2][Z3][Z4].

{\bf The Positive  Approach vs. the The Negative Approach}

We will present two {\it complementary} approaches to the enumeration of consecutive-Wilf classes,
both using the Umbral transfer matrix method.
The positive approach works better when you have many patterns, and the negative approach works better when there
are only a few, and works best when there is only one pattern to avoid.

{\bf Outline of the Positive Approach}

Instead of dealing with {\it avoidance} (the number of permutations that have zero occurrences of the given pattern(s))
we will deal with the more general problem of enumerating the number of permutations that have specified numbers
of occurrences  of {\it any} pattern of length $k$.

Fix a positive integer $k$, and let $\{ t_p: p \in S_k\}$ be $k!$ {\it commuting indeterminates} (alias variables).
Define the {\it weight} of an $n$-permutation $\pi=[\pi_1, \dots, \pi_n]$, to be denoted by $w(\pi)$, by:
$$
w([\pi_1, \dots , \pi_n]) :=
\prod_{i=1}^{n-k+1} t_{R([\pi_i,\pi_{i+1}, \dots, \pi_{i+k-1}])} \quad .
$$
For example, with $k=3$,
$$
w([2,5,1,4,6,3]) :=
t_{R([2,5,1])}t_{R([5,1,4])}t_{R([1,4,6])}t_{R([4,6,3])}=t_{231}t_{312}t_{123}t_{231}=t_{123}t_{231}^2t_{312} \quad .
$$
We are interested in an {\it efficient} algorithm for computing  the sequence of polynomials
in $k!$ variables
$$
P_n(t_{1 \dots k}, \dots, t_{k \dots 1 }):=
\sum_{\pi \in S_n} w(\pi) \quad ,
$$
or equivalently, as many terms as desired in the formal power series
$$
F_k(\{ t_p, p \in S_k\} ;z)=
\sum_{n=0}^{\infty} P_n z^n \quad .
$$

Note that once we have computed the $P_n$ (or $F_k$), we can answer {\it any} question about pattern avoidance
by specializing the $t$'s. For example  to get the number of $n$-permutations avoiding the single pattern
$p$, of length $k$, first compute $P_n$, and then plug-in $t_p=0$ and all the other t's to be 1.
If you want the number of $n$-permutations avoiding the set of patterns $\P$ (all of the same length $k$),
set $t_p=0$ for all $p \in \P$ and the other t's to be $1$. 
As we shall soon see, we will generate {\it functional equations} for $F_k$, featuring the
$\{ t_p \}$ and of course it would be much more efficient to specialize the $t_p$'s to 
the numerical values already in the functional equations, 
rather than crank-out the much more complicated
$P_n(\{ t_p \})$'s and then do the plugging-in.

First let's recall one of the many proofs that the number of $n$-permutations, let's denote it by $a(n)$, satisfies
the recurrence
$$
a(n+1)=(n+1)a(n) \quad .
$$
Given a typical member of $S_n$, let's call it $\pi=\pi_1 \dots \pi_n$, it can  be continued in $n+1$ ways, by
deciding on $\pi_{n+1}$. If $\pi_{n+1}=i$, then we have to ``make room'' for the new entry by incrementing
by $1$ all entries $\geq i$, and then append $i$. This gives a bijection between
$S_n \times [1,n+1]$ and $S_{n+1}$ and taking cardinalities yields the recurrence. Of course $a(0)=1$,
and ``solving'' this recurrence yields $a(n)=n!$. Of course this solving is ``cheating'', since
$n!$ is just shorthand for the solution of this recurrence subject to the initial condition $a(0)=1$, but
from now on it is considered ``closed form'' (just by convention!).

When we do {\it weighted counting} with respect to the weight $w$ with a given pattern-length $k$,
we have to keep track of the last $k-1$ entries of $\pi$:
$$
[\pi_{n-k+2} \dots \pi_{n}] \quad ,
$$
and when we append $\pi_{n+1}=i$, the new permutation
(let $a'=a$ if $a<i$ and $a'=a+1$ if $a \geq i$)
$$
\dots \pi_{n-k+2}' \dots \pi_{n}' i \quad ,
$$
has ``gained'' a factor  of $t_{R[\pi_{n-k+2}' \dots \pi_{n}' i]}$ to its weight.

This calls for the finite-state method, alas, the ``alphabet'' is indefinitely large, so we need
the umbral transfer-matrix method.

We introduce $k-1$ ``catalytic'' variables $x_1, x_2, \dots , x_{k-1}$, as well as a variable $z$
to keep track of the size of the permutation,
and $(k-1)!$ ``linear'' state variables $A[q]$ for each $q \in S_{k-1}$, to tell us
the state that the permutation is in.
Define the generalized weight $w'(\pi)$ of a permutation $\pi \in S_n$ to be:
$$
w'(\pi):=w(\pi) x_1^{j_1} x_2^{j_2} \dots x_{k-1}^{j_{k-1}} z^nA[q]\quad,
$$
where $[j_1, \dots , j_{k-1}]$, $(1 \leq j_1 < j_2 < \dots < j_{k-1} \leq n)$ is the {\it sorted}
list of the last $k-1$ entries of $\pi$, and $q$ is the reduction of its last $k-1$ entries.

For example, with $k=3$:
$$
w'([4,7,1,6,3,5,8,2])=
t_{231}t_{312}t_{132}t_{312}t_{123}t_{231}x_1^2x_2^8z^8A[21]=
t_{123}t_{132}t_{231}^2t_{312}^2x_1^2x_2^8 z^8 A[21]\quad .
$$

Let's illustrate the method with $k=3$. There are two states: $[1,2],[2,1]$ corresponding
to the cases where the two last entries are $j_1 j_2$ or $j_2 j_1$ respectively 
(we always assume $j_1<j_2$) .

Suppose we are in state $[1,2]$, so our permutation looks like

$$
\pi= [ \dots , j_1 , j_2 ]\quad ,
$$
and $w'(\pi)= w(\pi)x_1^{j_1} x_2^{j_2}z^nA[1,2]$.
We want to append $i$ ($1 \leq i \leq n+1)$ to the end. There are three 
cases.

{\bf Case 1}: $1 \leq i \leq j_1$ \quad .

The new permutation, let's call it $\sigma$, looks like
$$
\sigma= [ \dots j_1+1 , j_2+1 , i] \quad .
$$
Its state is $[2,1]$ and 
$w'(\sigma)=w(\pi)t_{231}x_1^{i}x_2^{j_2+1}z^{n+1}A[2,1]$.

{\bf Case 2}: $j_1+1 \leq i \leq j_2$ \quad .

The new permutation, let's call it $\sigma$, looks like
$$
\sigma= [ \dots j_1 , j_2+1 , i] \quad .
$$
Its state is also $[2,1]$ and $w'(\sigma)=w(\pi)t_{132}x_1^{i}x_2^{j_2+1}z^{n+1}A[2,1]$.

{\bf Case 3}: $j_2+1 \leq i \leq n+1$ \quad .

The new permutation, let's call it $\sigma$, looks like
$$
\sigma= [ \dots j_1 , j_2 , i] \quad .
$$
Its state is now $[1,2]$ and $w'(\sigma)=w(\pi)t_{123}x_1^{j_2}x_2^{i}z^{n+1}A[1,2]$.

It follows that any {\it individual} permutation of size $n$, and state $[1,2]$, gives rise to
$n+1$ children, and regarding weight, we have the ``umbral evolution'' (here $W$ is the fixed part of the
weight, that does not change):
$$
Wx_1^{j_1} x_2^{j_2} z^n A[1,2]
\rightarrow
Wt_{231} zA[2,1] \left ( \sum_{i=1}^{j_1}   x_1^i x_2^{j_2+1} \right ) z^n
$$
$$
+
W t_{132} zA[2,1] \left ( \sum_{i=j_1+1}^{j_2}   x_1^i x_2^{j_2+1} \right )  z^n
$$
$$
+
W t_{123} zA[1,2]  \left ( \sum_{i=j_2+1}^{n+1}   x_1^{j_2} x_2^{i} \right ) z^n \quad .
$$
Taking out whatever we can out of the $\sum$-signs, we have:
$$
Wx_1^{j_1} x_2^{j_2} z^n A[1,2]
\rightarrow
Wt_{231} zA[2,1] \left ( \sum_{i=1}^{j_1}   x_1^i  \right ) x_2^{j_2+1} z^n
$$
$$
+
W t_{132} zA[2,1] \left ( \sum_{i=j_1+1}^{j_2}   x_1^i  \right ) x_2^{j_2+1} z^n
$$
$$
+
W t_{123} zA[1,2]  \left ( \sum_{i=j_2+1}^{n+1}    x_2^{i} \right ) x_1^{j_2} z^n \quad .
$$

Now summing up the geometrical series, using the ancient formula:
$$
\sum_{i=a}^{b} Z^i =\frac{Z^a-Z^{b+1}}{1-Z} \quad ,
$$
we get
$$
Wx_1^{j_1} x_2^{j_2} z^n A[1,2]
\rightarrow
Wt_{231} zA[2,1] \left ( \frac{x_1-x_1^{j_1+1}}{1-x_1} \right ) x_2^{j_2+1} z^n
$$
$$
+
W t_{132} zA[2,1] \left ( \frac{x_1^{j_1+1}-x_1^{j_2+1}}{1-x_1} \right ) x_2^{j_2+1} z^n
$$
$$
+
W t_{123} zA[1,2]  \left ( \frac{x_2^{j_2+1}-x_2^{n+2}}{1-x_2} \right ) x_1^{j_2} z^n \quad .
$$
This is the same as:
$$
Wx_1^{j_1} x_2^{j_2} z^n A[1,2]
\rightarrow
Wt_{231} zA[2,1] \left ( \frac{x_1 x_2^{j_2+1} -x_1^{j_1+1}x_2^{j_2+1}}{1-x_1} \right ) z^n
$$
$$
+
W t_{132} zA[2,1] \left ( \frac{x_1^{j_1+1}x_2^{j_2+1}-x_1^{j_2+1}x_2^{j_2+1}}{1-x_1} \right ) z^n
$$
$$
+
W t_{123} zA[1,2]  \left ( \frac{x_1^{j_2} x_2^{j_2+1}-x_1^{j_2} x_2^{n+2}}{1-x_2} \right ) z^n \quad .
$$
This is what was called in [Z1], and its many sequels, a ``pre-umbra''. The above evolution
can be expressed for a general {\it monomial} $M(x_1,x_2,z)$ as:
$$
M(x_1,x_2,z)A[1,2]
\rightarrow
t_{231} zA[2,1] \left ( \frac{x_1 x_2 M(1,x_2,z) -x_1x_2M(x_1,x_2,z)}{1-x_1} \right ) 
$$
$$
+
 t_{132} zA[2,1] \left ( \frac{x_1x_2M(x_1,x_2,z)-x_1x_2M(1,x_1x_2,z)}{1-x_1} \right ) 
$$
$$
+
 t_{123} zA[1,2]  \left ( \frac{x_2M(1,x_1x_2,z)-x_2^{2}M(1,x_1, x_2z ) }{1-x_2} \right )  \quad .
$$
But, by {\it linearity}, this means that the coefficient of $A[1,2]$ (the weight-enumerator of all permutations of
state $[1,2]$) obeys the evolution equation:
$$
f_{12}(x_1,x_2,z)A[1,2]
\rightarrow
t_{231} zA[2,1] \left ( \frac{x_1 x_2 f_{12}(1,x_2,z) -x_1x_2f_{12}(x_1,x_2,z)}{1-x_1} \right ) 
$$
$$
+
 t_{132} zA[2,1] \left ( \frac{x_1x_2f_{12}(x_1,x_2,z)-x_1x_2f_{12}(1,x_1x_2,z)}{1-x_1} \right ) 
$$
$$
+
 t_{123} zA[1,2]  \left ( \frac{x_2f_{12}(1,x_1x_2,z)-x_2^{2}f_{12}(1,x_1,x_2z)}{1-x_2} \right )  \quad .
$$

Now we have to do it all over for a permutation in state $[2,1]$.
Suppose we are in state $[2,1]$, so our permutation looks like
$$
\pi= [ \dots , j_2 , j_1 ]\quad ,
$$
and $w'(\pi)= w(\pi)x_1^{j_1} x_2^{j_2}z^nA[2,1]$.
We want to append $i$ ($1 \leq i \leq n+1)$ to the end. There are three 
cases.

{\bf Case 1}: $1 \leq i \leq j_1$ \quad .

The new permutation, let's call it $\sigma$, looks like
$$
\sigma= [ \dots j_2+1 , j_1+1 , i] \quad .
$$
Its state is $[2,1]$ and 
$w'(\sigma)=w(\pi)t_{321}x_1^{i}x_2^{j_1+1}z^{n+1}A[2,1]$.

{\bf Case 2}: $j_1+1 \leq i \leq j_2$ \quad .

The new permutation, let's call it $\sigma$, looks like
$$
\sigma= [ \dots j_2+1 , j_1 , i] \quad .
$$
Its state is also $[1,2]$ and $w'(\sigma)=w(\pi)t_{312}x_1^{j_1}x_2^{i}z^{n+1}A[1,2]$.

{\bf Case 3}: $j_2+1 \leq i \leq n+1$ \quad .

The new permutation, let's call it $\sigma$, looks like
$$
\sigma= [ \dots j_2 , j_1 , i] \quad .
$$
Its state is now $[1,2]$ and $w'(\sigma)=w(\pi)t_{213}x_1^{j_1}x_2^{i}z^{n+1}A[1,2]$.

It follows that any {\it individual} permutation of size $n$, and state $[2,1]$, gives rise to
$n+1$ children, and regarding weight, we have the ``umbral evolution'' (here $W$ is the fixed part of the
weight, that does not change):
$$
Wx_1^{j_1} x_2^{j_2} z^n A[2,1]
\rightarrow
Wt_{321} zA[2,1] \left ( \sum_{i=1}^{j_1}   x_1^i x_2^{j_1+1} \right ) z^n
$$
$$
+
W t_{312} zA[1,2] \left ( \sum_{i=j_1+1}^{j_2}   x_1^{j_1} x_2^{i} \right )  z^n
$$
$$
+
W t_{213} zA[1,2]  \left ( \sum_{i=j_2+1}^{n+1}   x_1^{j_1} x_2^{i} \right ) z^n \quad .
$$

Taking out whatever we can out of the $\sum$-signs, we have:
$$
Wx_1^{j_1} x_2^{j_2} z^n A[2,1]
\rightarrow
Wt_{321} zA[2,1] \left ( \sum_{i=1}^{j_1}   x_1^i  \right )  x_2^{j_1+1} z^n
$$
$$
+
W t_{312} zA[1,2] \left ( \sum_{i=j_1+1}^{j_2}    x_2^{i} \right ) x_1^{j_1} z^n
$$
$$
+
W t_{213} zA[1,2]  \left ( \sum_{i=j_2+1}^{n+1}    x_2^{i} \right ) x_1^{j_1} z^n \quad .
$$
Now summing up the geometrical series, using the ancient formula:
$$
\sum_{i=a}^{b} Z^i =\frac{Z^a-Z^{b+1}}{1-Z} \quad ,
$$
we get
$$
Wx_1^{j_1} x_2^{j_2} z^n A[2,1]
\rightarrow
Wt_{321} zA[2,1] \left ( \frac{x_1-x_1^{j_1+1}}{1-x_1} \right ) x_2^{j_1+1} z^n
$$
$$
+
W t_{312} zA[1,2] \left ( \frac{x_2^{j_1+1}-x_2^{j_2+1}}{1-x_2} \right ) x_1^{j_1} z^n
$$
$$
+
W t_{213} zA[1,2]  \left ( \frac{x_2^{j_2+1}-x_2^{n+2}}{1-x_2} \right ) x_1^{j_1} z^n \quad .
$$
This is the same as:
$$
Wx_1^{j_1} x_2^{j_2} z^n A[2,1]
\rightarrow
Wt_{321} zA[2,1] \left ( \frac{x_1 x_2^{j_1+1} -x_1^{j_1+1}x_2^{j_1+1}}{1-x_1} \right ) z^n
$$
$$
+
W t_{312} zA[1,2] \left ( \frac{x_1^{j_1}x_2^{j_1+1}-x_1^{j_1}x_2^{j_2+1}}{1-x_2} \right ) z^n
$$
$$
+
W t_{213} zA[1,2]  \left ( \frac{x_1^{j_1} x_2^{j_2+1}-x_1^{j_1} x_2^{n+2}}{1-x_2} \right ) z^n \quad .
$$
The above evolution
can be expressed for a general {\it monomial} $M(x_1,x_2,z)$ as:
$$
M(x_1,x_2,z)A[2,1]
\rightarrow
t_{321} zA[2,1] \left ( \frac{x_1 x_2 M(x_2,1,z) -x_1x_2M(x_1x_2,1,z)}{1-x_1} \right ) 
$$
$$
+
 t_{312} zA[1,2] \left ( \frac{x_2M(x_1x_2,1,z)-x_2M(x_1,x_2,z)}{1-x_2} \right ) 
$$
$$
+
 t_{213} zA[1,2]  \left ( \frac{x_2M(x_1,x_2,z)-x_2^{2}M(x_1,1,x_2z)}{1-x_2} \right )  \quad .
$$
But, by {\it linearity}, this means that the coefficient of $A[2,1]$ (the weight-enumerator of all permutations of
state $[2,1]$) obeys the evolution equation:
$$
f_{21}(x_1,x_2,z)A[2,1]
\rightarrow
t_{321} zA[2,1] \left ( \frac{x_1 x_2 f_{21}(x_2,1,z) -x_1x_2f_{21}(x_1x_2,1,z)}{1-x_1} \right ) 
$$
$$
+
t_{312} zA[1,2] \left ( \frac{x_2f_{21}(x_1x_2,1,z)-x_2f_{21}(x_1,x_2,z)}{1-x_2} \right ) 
$$
$$
+
t_{213} zA[1,2]  \left ( \frac{x_2f_{21}(x_1,x_2,z)-x_2^{2}f_{21}(x_1,1,x_2z)}{1-x_2} \right )  \quad .
$$
Combining we have the ``evolution'':
$$
f_{12}(x_1,x_2,z)A[1,2]+ f_{21}(x_1,x_2,z)A[2,1] \rightarrow
$$
$$
t_{231} zA[2,1] \left ( \frac{x_1 x_2 f_{12}(1,x_2,z) -x_1x_2f_{12}(x_1,x_2,z)}{1-x_1} \right ) 
$$
$$
+ t_{132} zA[2,1] \left ( \frac{x_1x_2f_{12}(x_1,x_2,z)-x_1x_2f_{12}(1,x_1x_2,z)}{1-x_1} \right ) 
$$
$$
+
 t_{123} zA[1,2]  \left ( \frac{x_2f_{12}(1,x_1x_2,z)-x_2^{2}f_{12}(1,x_1,x_2z)}{1-x_2} \right )  \quad .
$$
$$
+t_{321} zA[2,1] \left ( \frac{x_1 x_2 f_{21}(x_2,1,z) -x_1x_2f_{21}(x_1x_2,1,z)}{1-x_1} \right ) 
$$
$$
+ t_{312} zA[1,2] \left ( \frac{x_2f_{21}(x_1x_2,1,z)-x_2f_{21}(x_1,x_2,z)}{1-x_2} \right ) 
$$
$$
+ t_{213} zA[1,2]  \left ( \frac{x_2f_{21}(x_1,x_2,z)-x_2^{2}f_{21}(x_1,1,x_2z)}{1-x_2} \right )  \quad .
$$
Now the ``evolved'' (new) $f_{12}(x_1,x_2,z)$ and $f_{21}(x_1,x_2,z)$ are the coefficients of $A[1,2]$ and $A[2,1]$
respectively, and since the {\it initial weight} of both of them is $x_1x_2^2z^2$, we have the
established the following system of functional equations:
$$
f_{12}(x_1,x_2,z)=x_1x_2^2z^2
$$
$$
+ t_{123} z \left ( \frac{x_2f_{12}(1,x_1x_2,z)-x_2^{2}f_{12}(1,x_1,x_2z)}{1-x_2} \right ) 
$$
$$
+ t_{312} z\left ( \frac{x_2f_{21}(x_1x_2,1,z)-x_2f_{21}(x_1,x_2,z)}{1-x_2} \right ) 
$$
$$
+ t_{213} z  \left ( \frac{x_2f_{21}(x_1,x_2,z)-x_2^{2}f_{21}(x_1,1,x_2z)}{1-x_2} \right )  \quad ,
$$
and
$$
f_{21}(x_1,x_2,z)=x_1x_2^2z^2
$$
$$
+ t_{231} z \left ( \frac{x_1 x_2 f_{12}(1,x_2,z) -x_1x_2f_{12}(x_1,x_2,z)}{1-x_1} \right ) 
$$
$$
+ t_{132} z \left ( \frac{x_1x_2f_{12}(x_1,x_2,z)-x_1x_2f_{12}(1,x_1x_2,z)}{1-x_1} \right ) 
$$
$$
+t_{321} z \left ( \frac{x_1 x_2 f_{21}(x_2,1,z) -x_1x_2f_{21}(x_1x_2,1,z)}{1-x_1} \right )  \quad .
$$

{\bf Let the computer do it!}

All the above was only done for {\it pedagogical} reasons. The computer can do it all automatically, much
faster and more reliably. Now if we want to find functional equations for the number of
permutations avoiding a given set of consecutive patterns $\P$, all we have to do is plug-in
$t_p=0$ for $p \in \P$ and $t_p=1$ for $p \not \in \P$. This gives a polynomial-time algorithm
for computing any desired number of terms. This is all done automatically in the
Maple package {\tt SERGI}. See the webpage of this article for lots of sample input and output.

Above we assumed that the members of the set $P$ are all of the same length, $k$.
Of course more general scenarios can be reduced to this case, where $k$ would be the largest length
that shows up in $P$.
Note that with this approach we end up with a
set of $(k-1)!$ functional equations in the $(k-1)!$ ``functions''
(or rather formal power series) $f_p$.

{\bf The Negative Approach }

Suppose that we want to quickly compute the first $100$ terms (or whatever) of the sequence enumerating
$n$-permutations avoiding the pattern $[1, 2, \dots, 20]$. 
As we have already noted,
using the ``positive'' approach, we have to set-up
a {\it system} 
of functional equations with $19!$ equations
and $19!$ unknowns. While the algorithm is still
{\it polynomial} in $n$ (and would give a ``Wilfian'' answer), it is not very practical! (This is yet another illustration why the
ruling paradigm in theoretical computer science, of equating ``polynomial time'' with ``fast'' is (sometimes) absurd).

This is analogous to computing words in a {\it finite} alphabet, say of $a$ letters,
avoiding a given word (or words) as
{\it factors} (consecutive subwords). If the word-to-avoid has length $k$, then the naive transfer-matrix
method would require  setting up a system of $a^{k-1}$ equations and  $a^{k-1}$ unknowns.
The elegant and powerful {\it Goulden-Jackson method} [GJ1][GJ2], beautifully exposited and extended in
[NZ], and even further extended in [KY], enables one to do it by solving one equation in
one unknown. We assume that the reader is familiar with it, and briefly describe the
analog for the present problem, where the alphabet is ``infinite''. This is also the approach
pursued in the beautiful human-generated papers [DK] and [KS]. We repeat that the {\it focus} and
{\it novelty} in the present work is in {\it automating} enumeration, and the current topic of
consecutive pattern-avoidance is used as a {\it case-study}.

First, some generalities!
For ease of exposition, let's focus on a single pattern $p$
(the case of several patterns is analogous, see [DK]).

Using the inclusion-exclusion ``negative'' philosophy for counting,
fix a pattern $p$. For any $n$-permutation, let $Patt_p(\pi)$ be the set
of occurrences of the pattern $p$ in $\pi$. For example
$$
Patt_{123} (179234568)=\{ 179, 234,345,456,568\} \quad ,
$$
$$
Patt_{231} (179234568)=\{ 792\} \quad ,
$$
$$
Patt_{312} (179234568)=\{ 923\} \quad ,
$$
$$
Patt_{132}(179234568)=Patt_{213}(179234568)=Patt_{321}(179234568)= \emptyset \quad .
$$

Consider the much larger set of pairs
$$
\{ (\pi, S) |\,\,  \pi \in S_n \quad \, \, , \, \, S \subset Patt_p(\pi) \},
$$
and define 
$$
weight_p(\pi,S):=(t-1)^{|S|} \quad,
$$
where $|S|$ is the number of elements of $S$.
For example,
$$
weight_{123} [179234568, \{ 234,568\}]=(t-1)^2 \quad ,
$$
$$
weight_{123} [179234568, \{ 179\}]=(t-1)^1=(t-1) \quad ,
$$
$$
weight_{123} [179234568, \emptyset]=(t-1)^0=1 \quad .
$$

Fix a (consecutive) pattern $p$ of length $k$, and
consider the weight-enumerator 
of all $n$-permutations according to the weight
$$
w(\pi):=t^{\# occurrences \,\, of \,\,  pattern  \,\, p \,\, in \,\, \pi} \quad ,
$$
let's call it $P_n(t)$. So:
$$
P_n(t):=\sum_{\pi \in S_n} t^{|Patt_{p}(\pi)|} \quad .
$$
Now we need the {\it crucial}, extremely deep, fact:
$$
t=(t-1)+1 \quad,
$$
and its corollary (for any finite set $S$):
$$
t^{|S|}=((t-1)+1)^{|S|}=\prod_{s \in S} ((t-1)+1)=
\sum_{T \subset S} (t-1)^{|T|} \quad .
$$
Putting this into the definition of $P_n(t)$, we get:
$$
P_n(t):=\sum_{\pi \in S_n} t^{|Patt_{p}(\pi)|} =
\sum_{\pi \in S_n} \sum_{T \subset Patt_p(\pi)} (t-1)^{|T|} \quad .
$$
This is the weight-enumerator (according to a different weight, namely $(t-1)^{|T|}$) of a much larger
set, namely the set of {\it pairs}, $(\pi, T)$, where $T$ is a subset of $Patt_p(\pi)$. Surprisingly, this is
much easier to handle!

Consider a typical such ``creature'' $(\pi,T)$. There are two cases

{\bf Case I}:
The last entry of $\pi$, $\pi_n$ does not belong to any of the members of $T$, in which case
chopping it produces a shorter such creature, in the set $\{1, 2, \dots , n \} \backslash \{ \pi_n \}$,
and reducing it to $\{1, \dots, n-1\}$ yields a typical member of size $n-1$. Since
there are $n$ choices for $\pi_n$, the weight-enumerator of creatures of this type
(where the last entry does not belong to any member of $T$) is $nP_{n-1}(t)$.

{\bf Case II}: Let's order the members of $T$ by their first (or last) index:
$$
[s_1,s_2, \dots, s_p] \quad ,
$$
where the last entry of $\pi$, $\pi_n$, belongs to $s_p$.
If $s_p$ and $s_{p-1}$ are disjoint, the ending cluster is simply $[s_p]$.
Otherwise $s_p$ intersects $s_{p-1}$. If $s_{p-1}$ and $s_{p-2}$ are disjoint, then
the ending cluster is $[s_{p-1},s_p]$. More generally,
the ending-cluster of the pair $[\pi,[s_1, \dots, s_p]]$
is the unique list $[s_i, \dots, s_p]$ that has the property that $s_i$ intersects $s_{i+1}$,
$s_{i+1}$ intersects $s_{i+2}$, $\dots$, $s_{p-1}$ intersects $s_p$, but
$s_{i-1}$ does not intersect $s_i$. It is possible that the ending-cluster of $[\pi,T]$ is the
whole $T$.

Let's give an example: with the pattern $123$.
The ending cluster of the pair:
$$
[157423689,[157,236,368,689]]
$$
is $[236,368,689]$ since $236$ overlaps with $368$ (in two entries) and $368$ overlaps with $689$ (also in two entries),
while $157$ is disjoint from $236$.

Now if you remove the ending cluster of $T$ from $T$ and remove the entries participating in the cluster from $\pi$,
you get a shorter creature $[\pi',T']$ where $\pi'$ is the permutation with all the entries in the
ending cluster removed, and $T'$ is what remains of $T$ after we removed 
that cluster. In the above example, we have
$$
[\pi',T']=[1574,[157]] \quad .
$$
Suppose that the length of $\pi'$ is $r$.

Let $C_{n}(t)$ be the weight-enumerator, according to the weight $(t-1)^{|T|}$, of canonical clusters
of length $n$, i.e. those
whose set of entries is
$\{1, \dots, n \}$. Then in Case II we have to choose a subset of $\{1, \dots, n \}$
of cardinality $r$ to be the $[\pi',T']$ and then choose a creature of size $r$ and a cluster of
size $n-r$. Combining Case I and Case II, we have , $P_{0}(t)=1$, and for $n \geq 1$:
$$
P_n(t)=nP_{n-1}(t)+ \sum_{r=2}^{n} {n \choose r} P_{n-r}(t)C_r(t) \quad .
$$
Now it is time to consider the {\it exponential generating function}
$$
F(z,t):=\sum_{n=0}^{\infty} \frac{P_n(t)}{n!}z^n \quad .
$$
We have
$$
F(z,t):=1+ \sum_{n=1}^{\infty} \frac{P_n(t)}{n!}z^n=
1+ \sum_{n=1}^{\infty} \frac{nP_{n-1}(t)}{n!}z^n+
\sum_{n=0}^{\infty} \frac{1}{n!} \left ( \sum_{r=2}^{n} {n \choose r} P_{n-r}(t)C_r(t) \right )z^n
$$
$$
=1+z \sum_{n=1}^{\infty} \frac{P_{n-1}(t)}{(n-1)!}z^{n-1}+
\sum_{n=0}^{\infty} \frac{1}{n!} \left ( \sum_{r=2}^{n} \frac{n!}{r!(n-r)!} P_{n-r}(t)C_r(t) \right )z^n
$$
$$
=
1+ z \sum_{n=0}^{\infty} \frac{P_{n}(t)}{n!}z^{n}+
\sum_{n=0}^{\infty} \left ( \sum_{r=2}^{n} \frac{1}{r!(n-r)!} P_{n-r}(t)C_r(t) \right )z^n
$$
$$
=1+ z F(z,t)+
\sum_{n=0}^{\infty} \left ( \sum_{r=2}^{n} \frac{P_{n-r}(t)}{(n-r)!} \frac{C_r(t)}{r!} \right )z^n
$$
$$
=1+ z F(z,t)+
\left ( \sum_{n-r=0}^{\infty} \frac{P_{n-r}(t)}{(n-r)!} z^{n-r} \right ) \left (  \sum_{r=0}^{\infty}  \frac{C_r(t)}{r!} z^r \right )
$$
$$
=1+ z F(z,t)+F(z,t)G(z,t) \quad,
$$
where $G(z,t)$ is the exponential generating function of $C_n(t)$:
$$
G(z,t):=\sum_{n=0}^{\infty} \frac{C_n(t)}{n!} z^n \quad .
$$
It follows that
$$
F(z,t)=1+ zF(z,t)+F(z,t)G(z,t) \quad,
$$
leading to
$$
F(z,t)=\frac{1}{1-z-G(z,t)} \quad .
$$
So if we would have a quick way to compute the sequence $C_n(t)$, we would have a quick way to compute 
the first {\it whatever} coefficients (in $z$) of $F(z,t)$ (i.e. as many $P_n(t)$ as desired).

\vfill\eject

{\bf A Fast Way to compute $C_n(t)$}

For the sake of pedagogy let the fixed pattern be $1324$. Consider a typical cluster
$$
[13254768,[1325,2547,4768]] \quad .
$$
If we remove the last atom of the cluster, we get the cluster
$$
[132547,[1325,2547]] \quad ,
$$
of the set $\{1,2,3,4,5,7\}$. Its canonical form, reduced to the set $\{1,2,3,4,5,6\}$, is:
$$
[132546,[1325,2546]] \quad .
$$
Because of the ``Markovian property'' (chopping the last atom of the clusters and reducing yields a shorter cluster),
we can build-up such a cluster, and in order to know how to add another atom, all we need to know
is the current last atom. If the pattern is of length $k$ (in this example, $k=4$), we need only to keep track of
the last $k$ entries. Let the sorted list (from small to large) be $i_1< \dots <i_k$, so the last
atom of the cluster (with $r$ atoms) is $s_r=[i_{p_1}, \dots, i_{p_k}]$,
where $1 \leq i_1 < i_2 <\dots < i_k \leq n$ is some increasing sequence of $k$ integers between $1$ and $n$.
We introduce $k$ {\it catalytic} variables
$x_1, \dots, x_k$, and define 
$$
Weight([s_1, \dots, s_{r-1}, [i_{p_1}, \dots, i_{p_k}]]):=z^n(t-1)^rx_1^{i_1} \cdots x_k^{i_k} \quad .
$$

Going back to the $1324$ example,
if we currently have a cluster with $r$ atoms, whose last atom is $[i_1, i_3, i_2, i_4]$, how can we add another atom?
Let's call it $[j_1, j_3, j_2, j_4]$
The new atom can overlap with the former one either in its last two entries, having:
$$
j_1=i_2 \quad j_3=i_4 \quad ,
$$
but because of the ``reduction'' (making room for the new entries) it is really
$$
j_1=i_2 \quad j_3=i_4+1 \quad ,
$$
(and $j_2$ and $j_4$ can be what they wish as long as $i_2<j_2<i_4+1<j_4 \leq n$).
The other possibility is that they only overlap at the last entry: 
$$
j_1=i_4 \quad
$$
(and $j_2, j_3, j_4 $ can be what they wish, provided that $i_4<j_2<j_3<j_4 \leq n$).

Hence we have the ``umbral-evolution'':
$$
z^n(t-1)^{r-1} x_1^{i_1}x_2^{i_2}x_3^{i_3}x_4^{i_4}
\rightarrow
z^{n+2}(t-1)^{r} \sum_{1 \leq j_1=i_2 < j_2 < j_3=i_4+1< j_4 \leq n} x_1^{j_1}x_2^{j_2}x_3^{j_3}x_4^{j_4}
$$
$$
+z^{n+3}(t-1)^{r} \sum_{1 \leq j_1=i_4 < j_2 < j_3< j_4 \leq n} x_1^{j_1}x_2^{j_2}x_3^{j_3}x_4^{j_4} \quad .
$$
These two iterated geometrical sums can be summed exactly, and from this ``pre-umbra'' the computer can deduce
(automatically!) the umbral operator, yielding
a functional equation for the {\bf ordinary} generating function
$$
\C(t,z; x_1, \dots , x_k)=\sum_{n=0}^{\infty} C_n(t;x_1, \dots, x_k) z^n \quad ,
$$
of the form
$$
\C(t,z; x_1, \dots , x_k)=(t-1)z^kx_1x_2^2 \dots x_k^k+
\sum_{\alpha} R_{\alpha}(x_1, \dots, x_k;t,z) \C(t,z; M_1^{\alpha}, \dots , M_k^{\alpha}) \quad,
$$
where $\{\alpha\}$ is a finite index set, $M_1^{\alpha}, \dots , M_k^{\alpha}$ are specific monomials
in $x_1, \dots, x_k,z$, derived by the algorithm, and $R_{\alpha}$ are certain rational functions
of their arguments, also derived by the algorithm.

Once again, the novelty here is that everything (except for the initial Maple programming) is done
{\it automatically} by the computer. It is the computer doing combinatorial research all on its own!

{\bf Post-Processing the Functional Equation}

At the end of the day we are only interested in $\C(t,z; 1, \dots , 1)$. 
Alas, plugging-in $x_1=1, x_2=1, \dots, x_k=1$ would give lots of $0/0$. Taking the limits,
and using L'H\^opital, is an option, but then we get a differential equation that would
introduce differentiations with respect to the catalytic variables, and we would not
gain anything.

But it so happens, in many cases, that the functional operator preserves some of the exponents of
the $x_i's$. For example for the pattern $321$ the last three entries are always $[3,2,1]$, and
one can do a {\it change of dependent variable}:
$$
\C(t,z; x_1, \dots , x_3)=x_1x_2^2x_3^3g(z;t)  \quad,
$$
and {\it now} plugging-in $x_1=1,x_2=1,x_3=1$ is harmless, and one gets 
a much simplified functional equation with {\it no} catalytic variables,
that turns out to be (according to S.B. Ekhad) the simple algebraic equation
$$
g(z,t)=-(t-1)z^2-(t-1)(z+z^2)g(z,t) \quad ,
$$
that in this case can be solved in closed-form (reproducing a result that goes back to [EN]).
Other times (like the pattern $231$), we only get rid of some of the catalytic variables.
Putting
$$
\C(t,z; x_1, \dots , x_3)=x_1x_2^2g(x_3,z;t)  \quad,
$$
(and then plugging-in $x_1=1,x_2=1$) gives a much simplified functional equation,
and now taking the limit $x_3 \rightarrow 1$ and using L'H\^opital (that Maple does all by itself)
one gets a pure differential equation for g(1,z;t), in $z$, that sometimes can be even
solved in closed form (automatically by Maple). But from the point of view of efficient
enumeration, it is just as well to leave it at that.

Any pattern $p$ is trivially equivalent to (up to) three other patterns (its reverse, its complement,
and the reverse-of-the-complement, some of which may coincide). It turns out
that out of these (up to) four options, there is one that is easiest to handle, and
the computer finds this one, by finding which ones gives the simplest functional (or if in luck
differential or algebraic) equation, and goes on to only handle  this representative.

{\bf The Maple package ELIZALDE}

All of this is implemented in the Maple package {\tt ELIZALDE}, that automatically produces {\it theorems}
and {\it proofs}. Lots of sample output (including computer-generated theorems and {\it proofs}) can
be found in the webpage of this article:

{\eighttt http://www.math.rutgers.edu/\~{}zeilberg/mamarim/mamarimhtml/auto.html} \quad .

In particular, to see all theorems and {\it proofs} for patterns of lengths 
$3$ through $5$ go to (respectively):

{\tt http://www.math.rutgers.edu/~zeilberg/tokhniot/sergi/oEP3\_200} ,

{\tt http://www.math.rutgers.edu/~zeilberg/tokhniot/sergi/oEP4\_60} ,

{\tt http://www.math.rutgers.edu/~zeilberg/tokhniot/sergi/oEP5\_40} .

If the proofs bore you, and by now you believe Shalosh B. Ekhad,
and you only want to see the statements of the {\it theorems}, for
lengths $3$ through $6$ go to (respectively):

{\tt http://www.math.rutgers.edu/~zeilberg/tokhniot/sergi/oET3\_200} ,

{\tt http://www.math.rutgers.edu/~zeilberg/tokhniot/sergi/oET4\_60} ,

{\tt http://www.math.rutgers.edu/~zeilberg/tokhniot/sergi/oET5\_40} .

{\tt http://www.math.rutgers.edu/~zeilberg/tokhniot/sergi/oET6\_30} .

Humans, with their short attention spans, would probably
soon get tired of even the statements of most of the theorems 
of this last file (for patterns of length $6$).

In addition to ``symbol crunching'' this package does quite a lot of ``number crunching''
(of course using the former). To see the ``hit parade'', ranked by size, together
with the conjectured asymptotic growth for single consecutive-pattern avoidance of lengths between $3$ and $6$,
see, respectively, the output files:

http://www.math.rutgers.edu/~zeilberg/tokhniot/sergi/oE3\_200 \quad ,

http://www.math.rutgers.edu/~zeilberg/tokhniot/sergi/oE4\_60 \quad ,

http://www.math.rutgers.edu/~zeilberg/tokhniot/sergi/oE5\_40 \quad ,

http://www.math.rutgers.edu/~zeilberg/tokhniot/sergi/oE6\_30 \quad .

Enjoy!

{\bf References}

[B] Andrew Baxter, {\it in preparation}.

[DK] Vladimir Dotsenko and Anton Khoroshkin, {\it Anick-type resolutions and consecutive pattern-avoidance},
arXiv:1002.2761v1[Math.CO] \quad .

[EKP] Richard Ehrenborg, Sergey Kitaev, and Peter Perry,
{\it A Spectral Approach to Consecutive Pattern-Avoiding Permutations},
arXiv: 1009.2119v1 [math.CO] 10 Sep 2010 \quad .

[EN] Sergi Elizalde and Marc Noy, {\it Consecutive patterns in permutations}, Advances in Applied Mathematics {\bf 30} (2003), 110-125.

[EZ] Shalosh B. Ekhad, and D. Zeilberger,
{\it Using Rota's Umbral Calculus to Enumerate Stanley's P-Partitions},
Advances in Applied Mathematics {\bf 41} (2008), 206-217.

[GJ1] Ian Goulden  and David M. Jackson,
{\it An inversion theorem for cluster decompositions of
sequences with distinguished subsequences},
J. London Math. Soc.(2){\bf 20} (1979), 567-576.
 
[GJ2] Ian Goulden  and David M. Jackson,
{\it ``Combinatorial Enumeration"}, John  Wiley, 1983, New York.

[KS] Anton Khoroshkin and  Boris Shapiro, {\it Using homological duality in consecutive pattern avoidance},
arXiv:1009.5308v1 [math.CO].

[KY]  Elizabeth J. Kupin and Debbie S. Yuster,
{\it Generalizations of the Goulden-Jackson Cluster Method}, J. Difference Eq. Appl. {\bf 16} (2010), 1563-5120.
arXiv:0810.5113v1[math.CO] .

[N] Brian Nakamura, {\it in preparation}.

[NZ] John Noonan and Doron Zeilberger,
{\it The Goulden-Jackson Cluster Method: Extensions, Applications, and Implementations},
J. Difference Eq. Appl. {\bf 5} (1999), 355-377.

[W] Herbert S. Wilf, {\it What is an answer}, Amer. Math. Monthly {\bf 89} (1982), 289-292.

[Z1] Doron Zeilbeger,
{\it The Umbral Transfer-Matrix Method I. Foundations},
J. Comb. Theory, Ser. A {\bf 91} (2000), 451-463 .

[Z2] Doron Zeilbeger,
{\it The Umbral Transfer-Matrix Method. III. Counting Animals}, 
 New York Journal of Mathematics {\bf 7}(2001), 223-231 .

[Z3] Doron Zeilbeger,
{\it The Umbral Transfer-Matrix Method V. The Goulden-Jackson Cluster Method for Infinitely Many Mistakes},
Integers {\bf 2} (2002), A5.

[Z4] Doron Zeilbeger,
{\it In How Many Ways Can You Reassemble Several Russian Dolls?}, Personal Journal of S.B. Ekhad and
D. Zeilberger, 
{\eighttt http://www.math.rutgers.edu/\~{}zeilberg/pj.html}, Sept. 16, 2009.

\end